\documentclass[twocolumn,10pt]{articlefederico}

\usepackage{amsthm}
\usepackage{amsmath}
\usepackage{amssymb}
\usepackage{color}

\usepackage{graphicx}

\newtheorem{theorem}{Theorem}

\newtheorem{axiom}[theorem]{Axiom}

\begin{document}
\title{\textsf{Todxs cuentan: cultivating diversity in combinatorics}}
\author{\textsf{Federico Ardila--Mantilla
\footnote{\noindent \textsf{San Francisco State University, San Francisco, CA,  USA, and Universidad de Los Andes, Bogot\'a, Colombia. federico@sfsu.edu \newline Supported by NSF CAREER Award DMS-0956178 and NIH SF BUILD grant 1UL1MD009608-01.}}}}

\date{}
\maketitle

After nine years of activity, the SFSU-Colombia Combinatorics Initiative has helped build a strong, active community of more than 200 mathematicians, most of whom are members of underrepresented groups in mathematics. 
More than 50 of them have pursued Ph.D.s in mathematics, while others continue to be mathematics users, enthusiasts, and ambassadors in other fields, and to encourage and inspire the next generation of scientists in their communities. 
This article tells our story, and shares some lessons we have learned about broadening and deepening representation in mathematics.

\begin{figure}[h]
\begin{center}
\includegraphics[width=3in]{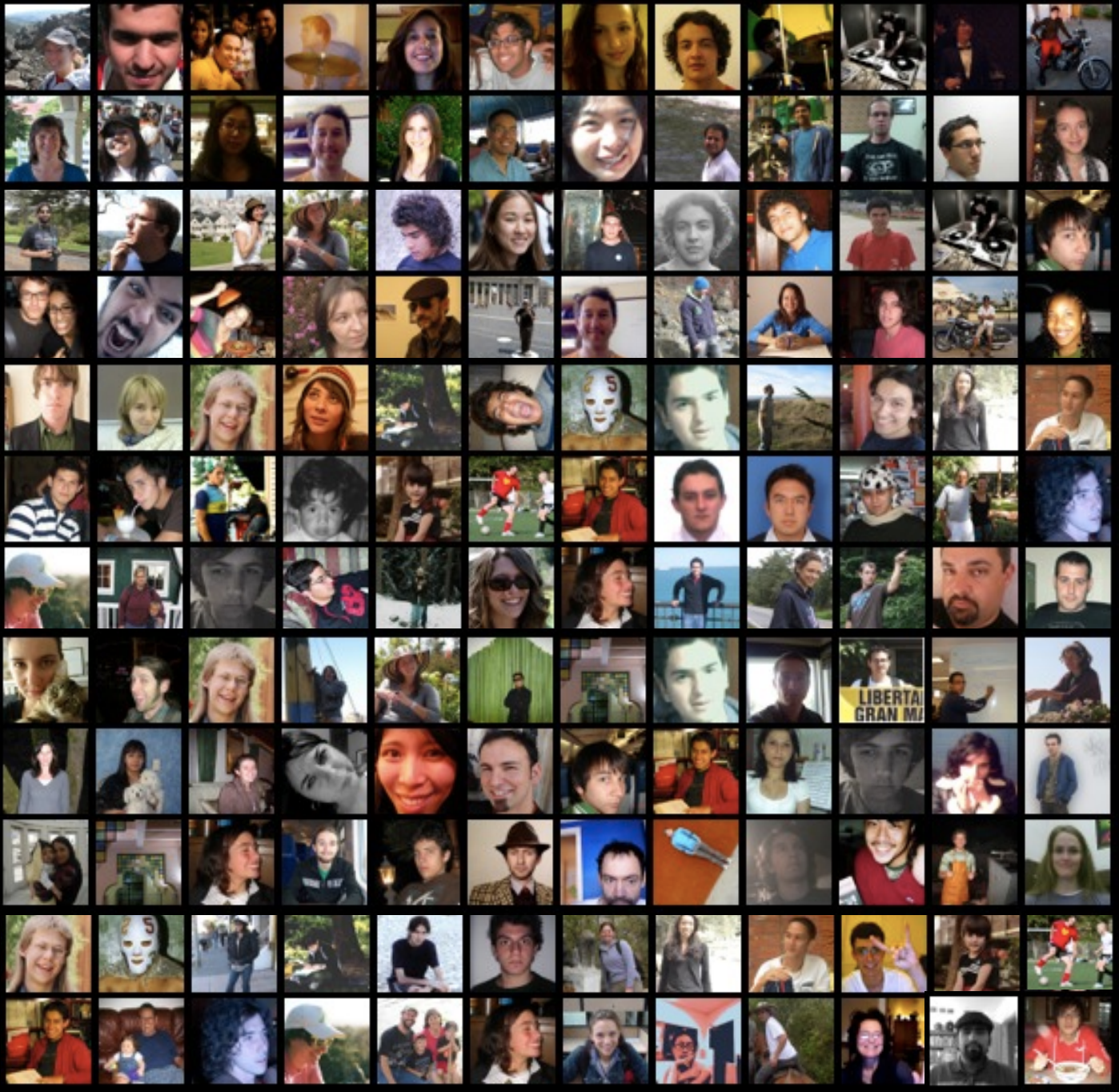} 
\caption{The SFSU-Colombia combinatorics family.}
\end{center}
\end{figure}

\vspace{-.7cm}
%

We begin by stating our axioms in Section 0. We then outline our work in Section 1, and discuss our underlying sociopolitical framework and pedagogical strategies in Section 2.



\smallskip
\noindent
\begin{small}
\textbf{\textsf{0. THE AXIOMS}}
\end{small}
\smallskip

%
Let me begin with some axioms. I cannot prove them but I firmly believe in them, and I build my work upon them.

\begin{axiom}
Mathematical potential is distributed equally among different groups, irrespective of geographic, demographic, and economic boundaries.
\end{axiom}

\begin{axiom}
Everyone can have joyful, meaningful, 
and empowering mathematical experiences.
\end{axiom}

\begin{axiom}
Mathematics is a powerful, malleable tool that can be shaped and used differently by various communities to serve their needs.
\end{axiom}

\begin{axiom}
Every student deserves to be treated with dignity and respect.
\end{axiom}

These statements should not sound revolutionary; but considering the current practices of the mathematical society, they are a call to action.

\medskip
\noindent
\begin{small}
\textbf{\textsf{1. SFSU-COLOMBIA COMBINATORICS INITIATIVE}} 
\end{small}

\smallskip

The SFSU-Colombia Combinatorics Initiative is a research and training collaboration  which seeks to offer every interested student a challenging and supportive mathematical experience, while planting seeds for the broader and deeper representation of different groups in mathematics.


\smallskip
\noindent
\textbf{\textsf{The partners.}} 
Los Andes is an elite private university. Its math department is one of Colombia's strongest, and through scholarships, it attracts some of the best prepared students in the country. Most Los Andes participants in the program are undergraduates with a very solid mathematics background but little understanding of what research looks like. About one fifth of them were members of the Colombian math olympiad team, which I coached for 15 years.

SFSU is a large public 4-year university with a diverse population: over 60\% of the students come from ethnic minority groups, and almost half are first-generation college students. The Master's program welcomes and serves a student body with a wide range range of demographics and academic preparation. 
SFSU is home to an active research group in combinatorics, including about 15 students at a time.
 
I was born and raised in Colombia, discovering mathematics through the Olimpiadas Colombianas de Matem\'aticas. I came to the US 22 years ago as an undergraduate on a scholarship at MIT and have been here since, while remaining in close contact with Colombia and its mathematics. In the US I am usually counted as a minority mathematician, and I have often felt in the minority. 
Though I was treated well, most of my training in the US took place alone, avoiding mathematical spaces where I felt uncomfortable.
However, unlike most students from marginalized groups in mathematics, I never had to overcome the structural inequalities of our educational system. In particular, I never had a teacher or peer who doubted my abilities or told me that I was not good enough to succeed. 

\smallskip
\noindent
\textbf{\textsf{How it all started (2005)}}.
Felipe Rinc\'on, then an undergraduate at Los Andes, had taken my Algebraic Combinatorics course there in 2003, and was writing his thesis under my supervision. In view of his classmates' interest and the relative lack of activity in this field in Colombia, he volunteered to teach an unofficial combinatorics course for free. The course was a great success, attracting more than 20 of the approximately 120  math majors.

I had recently begun an Assistant Professorship at San Francisco State University, and was very interested in contributing to mathematics in Colombia. Felipe's course exposed a great need, so I decided to offer a followup topics course on matroid theory, offered jointly at SFSU in person and at Los Andes electronically. Internet courses were still not in fashion, and this was exactly the kind of wild, uncertain experiment that Colombians love to embark upon. 

\smallskip

\noindent
\textsf{{\textbf{SFSU/Los Andes, Matroid Theory (2007)}}}.
%
%
%
%
I taught the class at SFSU and filmed it using my colleague Arek Goetz's  wonderfully low-budget artisanal setup, modeled after the one used by his figure skating sister: heat sensors across the front of the room detected where I was standing, and told the camera in the back where to point. With practice I learned to use large and clear handwriting, and to move around less so I would not make the viewers dizzy. Arek also improved his algorithm, so the camera would not follow me when I walked across the room and back to get the eraser.

Los Andes students watched the lecture together in a classroom, where I was present virtually to answer questions. We could not meet simultaneously due to time zone differences and the limits of the technology available to us. I visited Bogot\'a early in the semester to meet the students in person. Ph.D. students from UC Berkeley also took the class;
some in person, and some on video.

We made great efforts for  US and Colombian students to feel that they were in the same class. Students created minibios including their photos, personal background, and mathematical interests; they were not afraid to let their personality shine through. I also created an online forum where US and Colombian, got to `meet',  
discuss the course material and the assignments, and find future collaborators. 
\begin{figure}[h]
\includegraphics[width=3.15in]{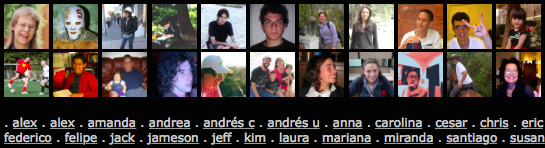}
\caption{Student profiles from the 2007 joint course.}
\end{figure}


In the first half of the course I assigned homework which ranged from reasonably straighforward exercises, accessible to anyone in the class, to approachable but  challenging recent results in the matroid theory literature, which required extensive (sometimes international) teamwork. Occasionally I assigned unannounced open problems. Posting student solutions -- while trying to represent everyone in the class -- raised the quality of the work and the writing. 
Most but not all Colombians were proficient in English, and some took the opportunity to write mathematics in English for the first time. Others didn't, and SFSU students got to read mathematics in Spanish; this was especially exciting for some of the US Latinas/os, many of whom spoke Spanish at home but not in their mathematical life.

In the second half, students did final projects in pairs; I suggested dozens of projects, from surveys of classic topics to current open problems of interest. Most students tried to do original research.

\begin{figure}[h]
\includegraphics[height=1in]{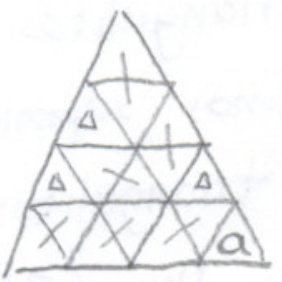} \quad
\includegraphics[height=1.2in]{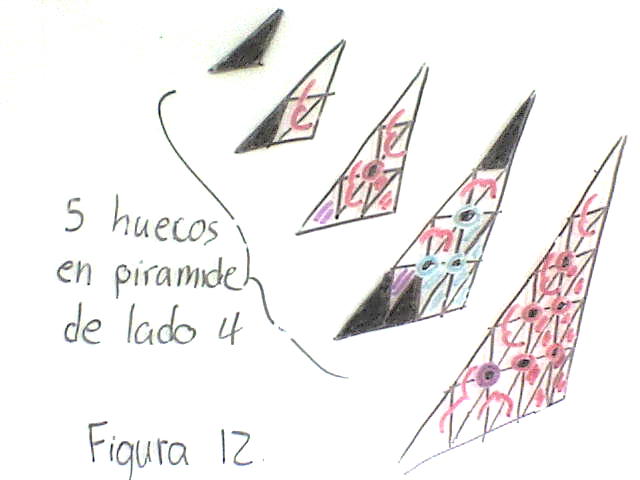} 
\includegraphics[height=1.7in]{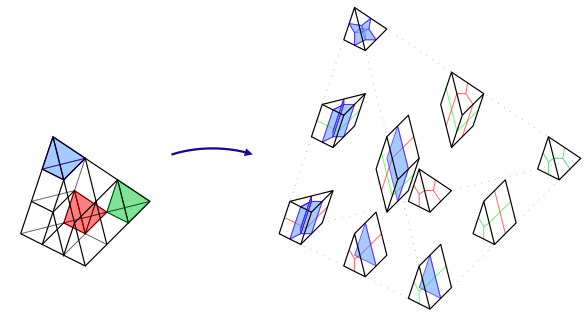}

\caption{Cesar Ceballos's work progresses from a homework assignment to a final project, a Master's thesis, and a publication \cite{AC} about the triangulations of the product of two simplices $\Delta_m \times \Delta_n$.}
\end{figure}

\vspace{-.3cm}

\noindent 
$\bullet$
\emph{Results.} 
The 11 projects in the class led to 5 papers in international journals (3 coauthored internationally) \cite{ABD, AC, AFR, AR, DF} and 8 theses. 
Nine years later, of the 21 students in the class, 11 have completed Ph.D.s in mathematics. Currently,


\noindent
-- 4 are university professors in mathematics,


\noindent
-- 6 are postdoctoral researchers in mathematics,



\noindent
-- 4 are community college faculty in mathematics,


\noindent
-- 2 are Ph.D. students in economics, and


\noindent
-- 5 are working in industry.


\smallskip

\noindent 
$\bullet$
\emph{Observations and lessons learned.}

\noindent
-- Every student did deep work, especially in their final projects. Students with substantial mathematical gaps learned what they needed along the way.



\noindent 
-- A strong sense of teamwork and collaboration helped most students lift each other up, though a few struggled with having to do mathematics in groups.


\noindent
-- Most projects were international, and those were most productive. 
Students' curiosity towards and accountability to strangers played a useful role.

\noindent
-- Generally speaking, SFSU students were impressed with the knowledge and problem-solving skills of Los Andes students. Los Andes students were impressed by the work ethic and determination of SFSU students, particularly in their research projects.

\noindent
-- Asynchronicity was far from ideal; it made class lecture-centric. Creating a coherent shared experience took hard work, and it was still not a great substitute for personal interaction.



%
%
%
%
%

\begin{figure}[h]
\includegraphics[width=3.15in]{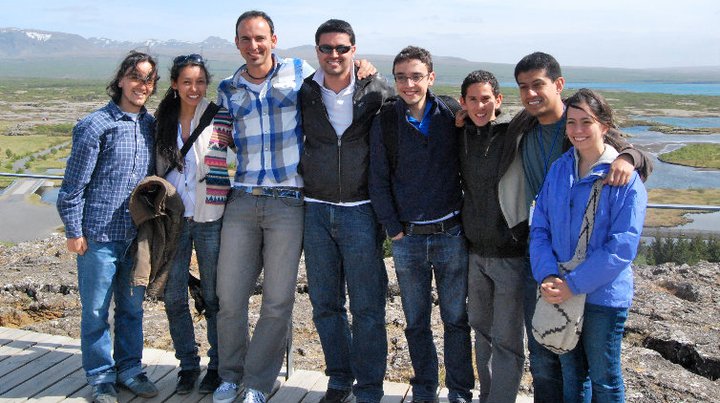} 
\caption{With students presenting their work at the FPSAC 2011 combinatorics conference in Iceland.}
\end{figure}

\smallskip

\noindent
\textbf{\textsf{Some guiding principles for the course design.}} 
The Matroid Theory class became the blueprint for 8 topics courses over the last 9 years; 6 of them were offered jointly at SFSU and Los Andes.
As I continued to design these courses, I identified a few key principles that I always attempted to implement.

\noindent
 --  Choose a deep topic of current interest that is accessible to students with different backgrounds.

\noindent
 --  Hold students to extremely high standards, and match that expectation with a  solid support system.

\noindent
 --  Devise challenging, interesting, and inspiring assignments, including a final project in pairs that students have the freedom to design themselves.

\noindent 
 --  Give students the time they need: allow 2 weeks for each homework and 2 months for the final project.

%


\noindent
 -- 
Create a course structure that builds a strong community through a shared mathematical purpose.

\noindent Several student evaluations mentioned this sense of community, and some took it to the social media:

\begin{figure}[h]
\includegraphics[width=3.15in]{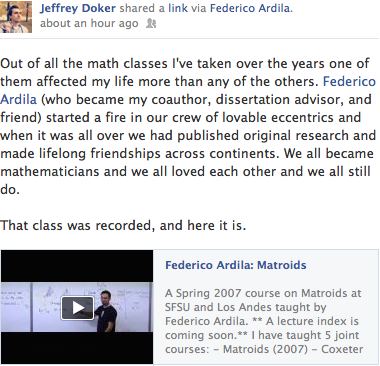} 
\end{figure}

\smallskip
\noindent
\textbf{\textsf{Online resources.}} I made all materials for my 6 SFSU-Colombia courses freely available, including more than 240 hours of videos on the YouTube Channel \texttt{federicoelmatematico}. According to the channel's statistics, users in 155 countries have viewed more than 10,000 hours of combinatorics videos. More meaningfully to me, I have personally heard from active users in places like Colombia, Germany, India, Iran, Sudan, Turkey, and the US.

Technology offers exciting possibilities to increase access to education, but we should be cautious in interpreting the numbers above. For example, the average view duration for these 50-75 minute videos was 8:33. Of the viewers who divulged their gender, 81\% were men. The last few lectures of each course were watched between 60 and 300 times each, with viewers concentrated in the US and Western Europe, and including several of my students and colleagues. These online resources have certainly been helpful to people, but I do not believe they have had a significant effect in closing the gap for access to quality mathematics resources.

%
%
%
%
%
%
%
%
%
%
%

%
%
%

\begin{figure}[h]
\begin{center}
\includegraphics[width=3.15in]{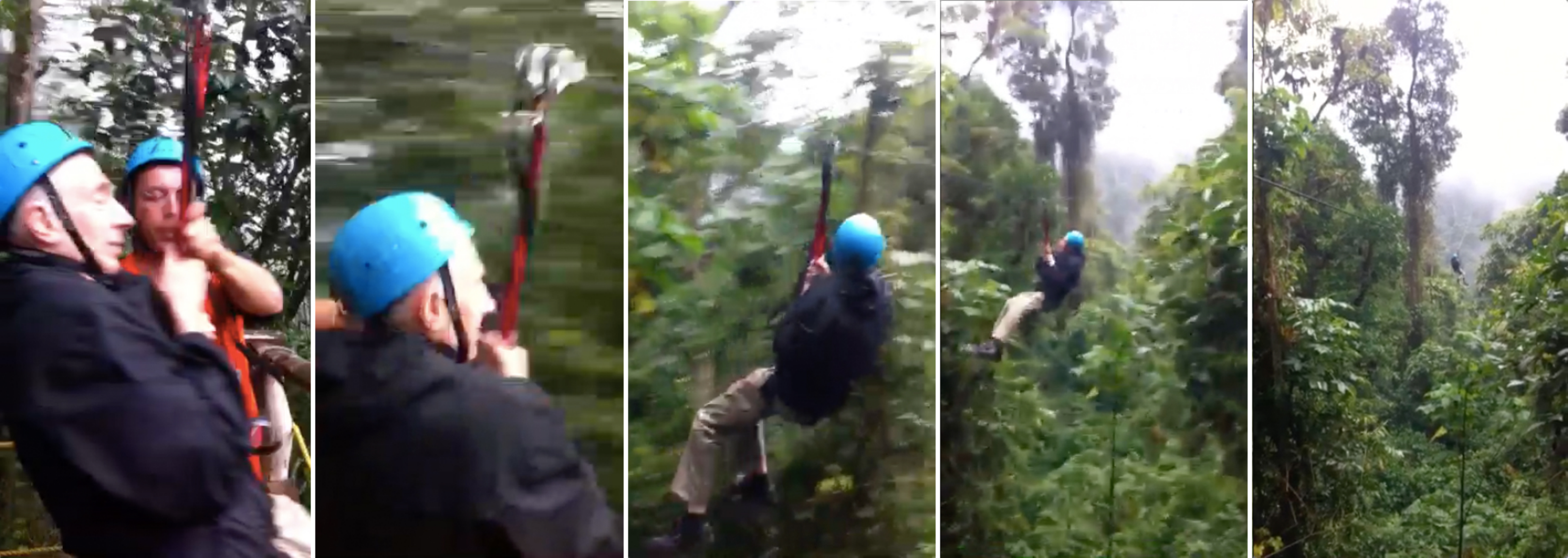} 
\caption{MIT professor Richard Stanley, as he is about to zipline across the Chicaque Cloud Forest with US and Colombian students: ``I hope I don't have a heart attack right now;
I still have to go to my Birthday Conference next week!"}
\end{center}
\end{figure}
%
%
%
\vspace{-0.5cm}

\noindent
\textbf{\textsf{Encuentro Colombiano de Combinatoria (ECCO).}} 
It soon became clear that US and Colombian students wanted to build closer ties and collaborate in person. 
I organized the Encuentro Colombiano de Combinatoria to benefit young mathematicians first and foremost, especially those who do not have easy access to such an opportunity. Now organized by a committee of former participants, ECCO meets biyearly. We do our best to build an atmosphere that is equally professional and welcoming.

ECCO features minicourses by international experts, collaborative problem sessions, research talks by students, open problem presentations, mentoring sessions, a hike or two, and (inevitably, it seems) an impromptu street party.


Now that ECCO has gained some international notoriety, several combinatorics experts have asked to attend. We welcome them with care, keeping in mind that this is a school and an \emph{encuentro}\footnote{gathering, usually of people with shared experiences}, not a regular conference. We ask these experts to do problem sets with the students, to present open problems that they would like help with, to serve as mentors, and probably to join the dance floor at some point.


I consider it a success that I am no longer on the ECCO organizing committee; 
Sharing the decision making has brought new energy, ideas, and perspectives, and has given many people a sense of ownership of this project. More importantly, I believe it has fostered their agency as members of the mathematical community, and empowered them to pursue their own initiatives.



\begin{figure}[h]
\includegraphics[width=3.15in]{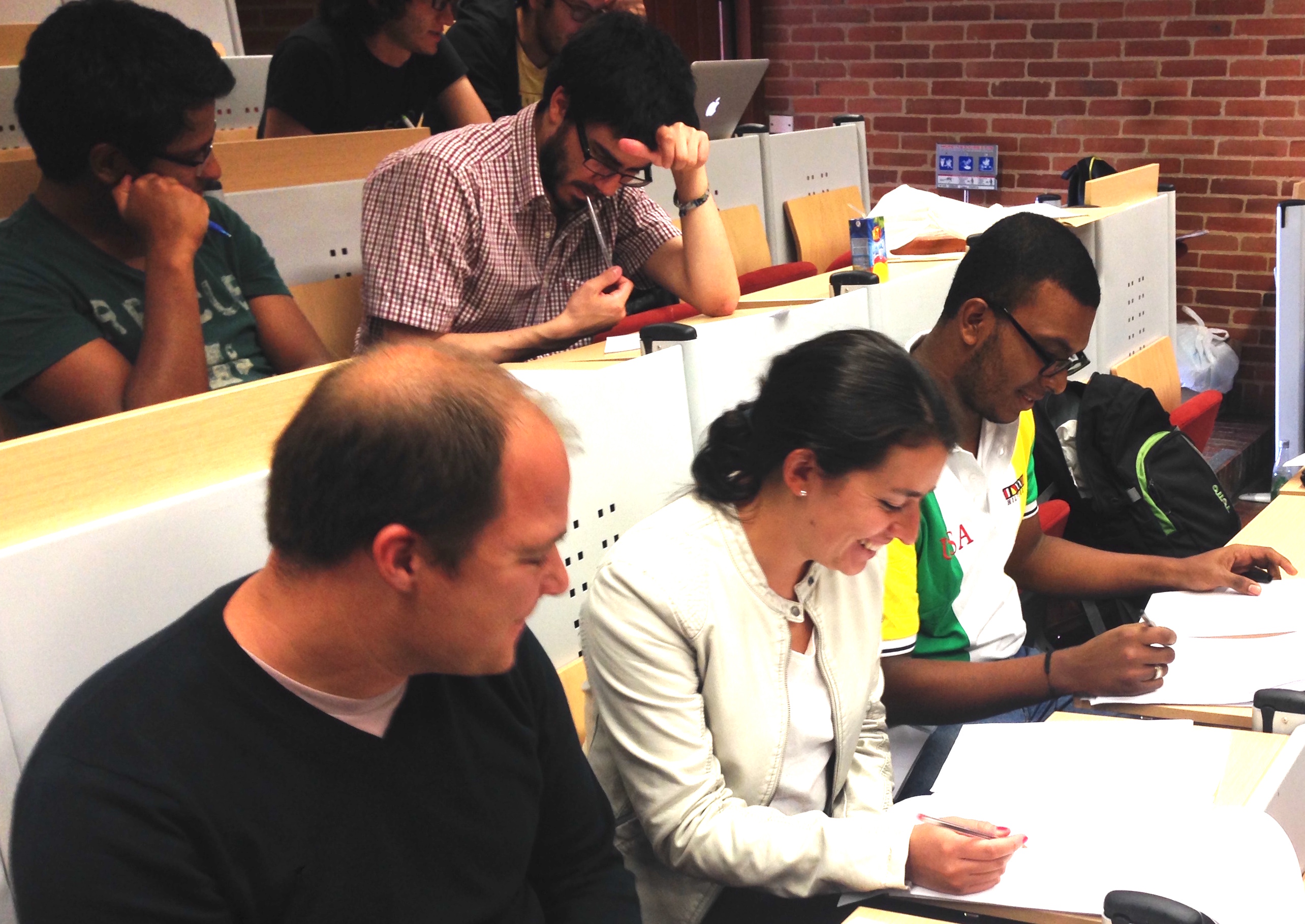} 
\caption{Andr\'es Salda\~na (U. del Valle) learning combinatorics and speaking English, both for the first time, with Brian Davis (SFSU) and Eliana Duarte (Los Andes) at ECCO 2014. 
He is now a graduate student at the Berlin Mathematical School.} 
\end{figure}

\noindent
\textbf{\textsf{Funding.}} The Matroid Theory class and ECCO 2008 were generously funded by seed grants from SFSU and Los Andes. With those results in hand, I applied for the NSF CAREER and RUI grants that have funded our activities since then.

The cost of tuition and living in San Francisco is a real challenge for our students, most of whom have significant work loads while they study. We have supported some of them using other NSF and NIH grants at SFSU, such as Matt Beck's GK-12 grant in mathematics and Frank Bayliss's SEO grants in science. 
We plan to implement a more sustainable funding infrastructure that will allow our students to better focus their efforts on their academic preparation.

\smallskip
\noindent
\textbf{\textsf{Some statistics.}}
Approximately 200 students have officially enrolled in the classes offered by this initiative, and approximately 50 of them have pursued Ph.D.s in the mathematical sciences. 

The initiative has had 40 theses students (4 Ph.D., 27 M.S., 12 B.S.): 28 in the US (14 women, 15 URMs\footnote{underrepresented minorities}) and 12 in Colombia (3 women, 12 URMs), who authored more than 20 publications, including journals such as \emph{Duke Math. J., Advances in Math, Int. Math. Res. Not., and J. Combin. Th. (A).} 

Of these 40 students, 28 have entered math Ph.D.s (19 US, 11 women, 21 URMs), and 25 have finished or are current students, including all 21 URMs.\footnote{Numbers for US and women omitted for confidentiality.}



\medskip
\noindent
\begin{small}\textbf{\textsf{2. DEEPENING REPRESENTATION.}} 
\end{small}
\smallskip

The underrepresentation of women, Latinas/os, and African-Americans  in US mathematics is well documented.\footnote{We focus for the moment on these three large groups for which data is available.} These groups constituted respectively

\noindent $\bullet$
50.8\%, 17.4\%, 13.2\%  of the 2014 population \cite{UScensus},

\noindent $\bullet$
31\%, 3.5\%, 2.5\% of the 2014 new math Ph.D.s \cite{AMSSurvey},

\noindent $\bullet$
18\%, 3\%, 1\% of the full-time math faculty at Ph.D. granting institutions in 2005. \cite{CBMSSurvey} 

Underrepresentation is more drastic further along the academic career.
Among the science, engineering, and health faculty in 2008, Latinos and African-Americans constituted 4.8\% and 5\% of assistant professors, 
3.6\% and 4.8\% of associate professors, 
and 3.3\% and 2.3\% of full professors. \cite{NSFreport} 

These numbers naturally lead to further underrepresentation in positions of leadership and decision making power, 
slowing down the changes necessary to reverse this trend.
Within the AMS, all 63 presidents have been white and 61 have been men. The situation is improving: the 12 current officers include 4 women and 2 foreign-born Latinos. A new Directorship of Education and Diversity was instituted in 2016, and at a more grassroots level, many excellent Mathematics Programs That Make a Difference have been recognized since 2006. However, the mathematical society at large still has a lot of work to do.

%
%
%
%
%

It appears that, now more than ever, the US society recognizes the importance of building a scientific workforce that reflects and represents our diverse demographics. It is less clear how this will be achieved, but we have learned a few things: it is not enough to wait for the occasional lone `geniuses' from unlikely places to make it against all odds; it is not enough to offer a few different kinds of people entrance to the same old house, and expect them to come in and to thrive. To truly broaden representation we must deepen representation. We must be prepared to truly accept the structural inequalities that led to the current state of affairs, both within our mathematical community and in our nation at large, and to do the deep, hard work that is required to overcome them. As a mathematician with rather limited knowledge and expertise, I find this a daunting challenge. 
%
%

As scientists do when faced with a seemingly intractable problem, I have focused  on a smaller problem for the moment: building nurturing environments where different people can thrive and do excellent mathematics. I certainly don't have a recipe for how to do this, and I know we still have a long way to go, but I do think the SFSU-Colombia Combinatorics Initiative has achieved some success.
Along the way, by listening to students, teachers, organizers, and scholars, I have learned a few things about the larger structural obstacles to be overcome. I believe the following factors have been crucial to the small successes we've had so far.

%


\smallskip

\noindent 
\textsf{\textbf{The mathematics.}} To succeed in mathematics, one must do high quality mathematics. This is especially true for students without elite credentials, for whom the bar to success is often set higher. It is crucial to involve mentors with active research programs in areas of current interest.

\smallskip

\noindent 
\textsf{\textbf{Todos cuentan.}} 
Our activities are designed to serve \textbf{every} interested student, by building an inclusive environment that everyone contributes to and benefits from. We aim to increase students' sense of belonging and responsibility to their mathematical community.

Intentionally or not, most programs in higher mathematics are designed to select the `top' mathematicians at each stage, and prepare them for the next stage. Such programs have certainly played a crucial role in the careers of many mathematicians, including mine. However, they have also played a role in excluding and discouraging others with great mathematical potential, particularly among underrepresented groups. 


\smallskip

\noindent 
\textsf{\textbf{Equitable spaces.}}
We try to be mindful of how different students experience the same environment, and find ways to make sure each one of them is actively involved and engaged within the spaces we provide. 

Language matters. Many of the standard patterns of communication of mathematicians -- like calling a fact `easy to see' when it isn't, or saying 'I get it' when one doesn't -- are harmless to some, but they feel exclusive to those who already feel like outsiders. 

Structure matters. Without an explicit and mindful effort, classrooms easily turn into  conversations between a professor and a handful of students. My current courses are organized so that every student participates in every class, either verbally, in writing, or in small groups; I use several of Kimberly Tanner's techniques for creating equitable learning environments.  \cite{T}


\begin{figure}[h]
\includegraphics[width=3.7in]{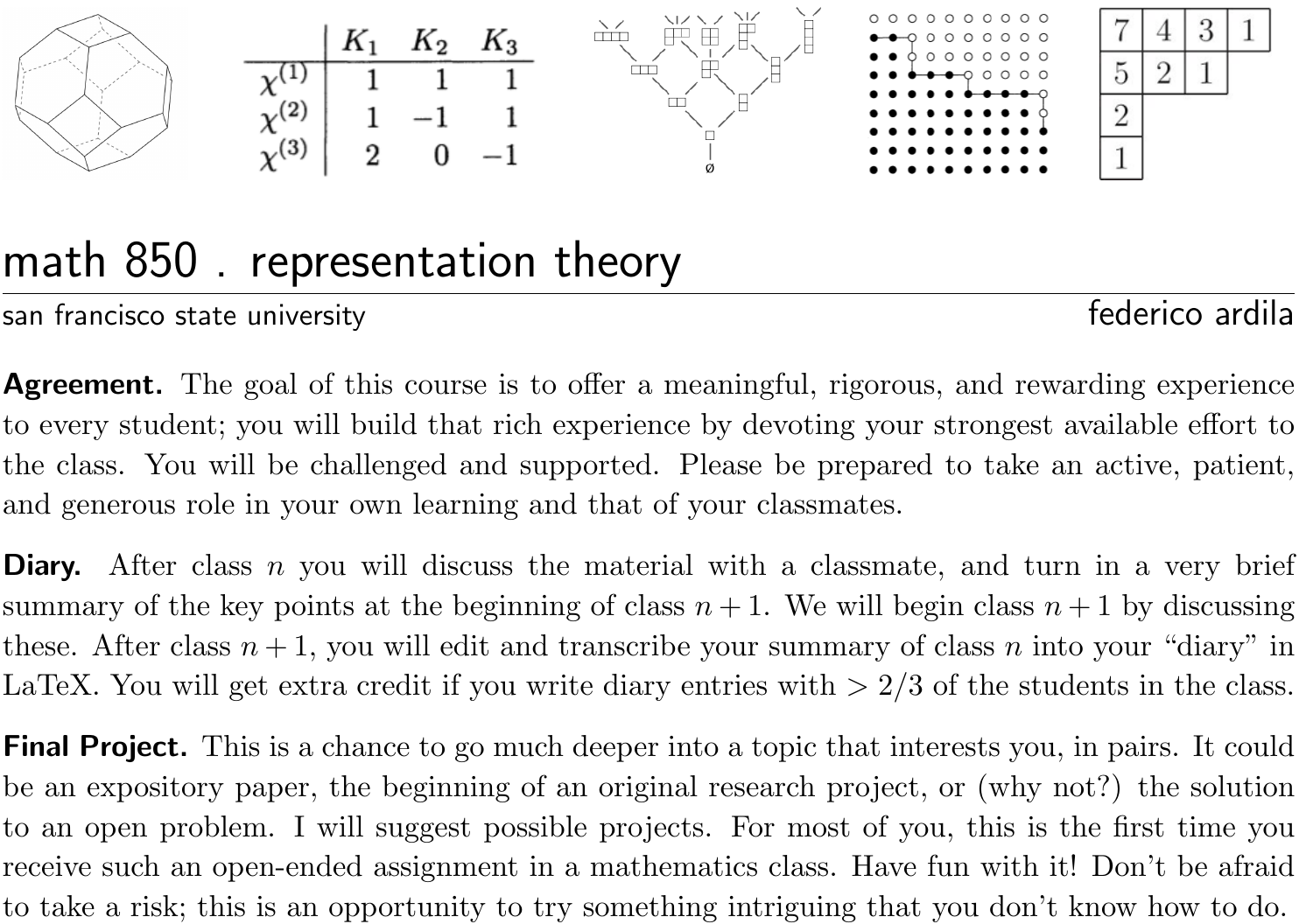} 
\caption{Language and structure: syllabus excerpts.}
\end{figure}

%

\smallskip

\noindent 
\textsf{\textbf{The support system.}} 
A challenging experience 
can easily become alienating if it is not presented mindfully and accompanied with abundant support. 
We do not shy away from the emotional side of this work. A career in mathematics requires balancing long periods of  frustration with (sometimes too brief) moments of great joy, and students find it surprising and beneficial to learn that they are not the only ones struggling with that balance. Readings and discussions on the psychology of mathematics and science \cite{Sc,Th}, 
growth mindset \cite{D}, stereotype threat \cite{S}, and impostor syndrome \cite{CI} have been helpful. 

At SFSU we collaborate with initiatives like SF BUILD, which works to create inclusive and supportive research environments across 6 departments, and the Mathematistas student group for women in math, which builds community in our department through many informal activities. 
 Outside of SFSU we benefit from interacting with organizations such as SACNAS, the Math Alliance, USTARS, and Latin@s in Math.


We are aware that the crucial work of supporting students traditionally falls  mostly on women and people of color \cite{J, MLHA, MOWA}. We make every attempt to counter that tendency. 

%
%
%

\smallskip

\noindent 
\textsf{\textbf{People rise to high expectations.}}
Many of my students from marginalized groups have been told, often by well-meaning professors in positions of power over them, that they cannot do something or that they are not good enough to be mathematicians. I \textbf{never} say this to a student. I cannot possibly know that.\footnote{Upon meeting Endre Szemer\'edi, 
Israel Gelfand told him ``Just try to find another profession; there are plenty in the world where you may be successful." \cite{Sz} Szemer\'edi went on to 
write his Ph.D. thesis under Gelfand, and 
win the 2012 Abel Prize among many other honors.}

I have worked with students whose mathematical potential is not immediately apparent to me, but I know they are here for a reason. 
My approach is to always treat them with respect, give them an intriguing project that suits their experience and interests, and support them along the way; I have seen most of them rise to that challenge.

\begin{figure}[h]
\includegraphics[width=3.15in]{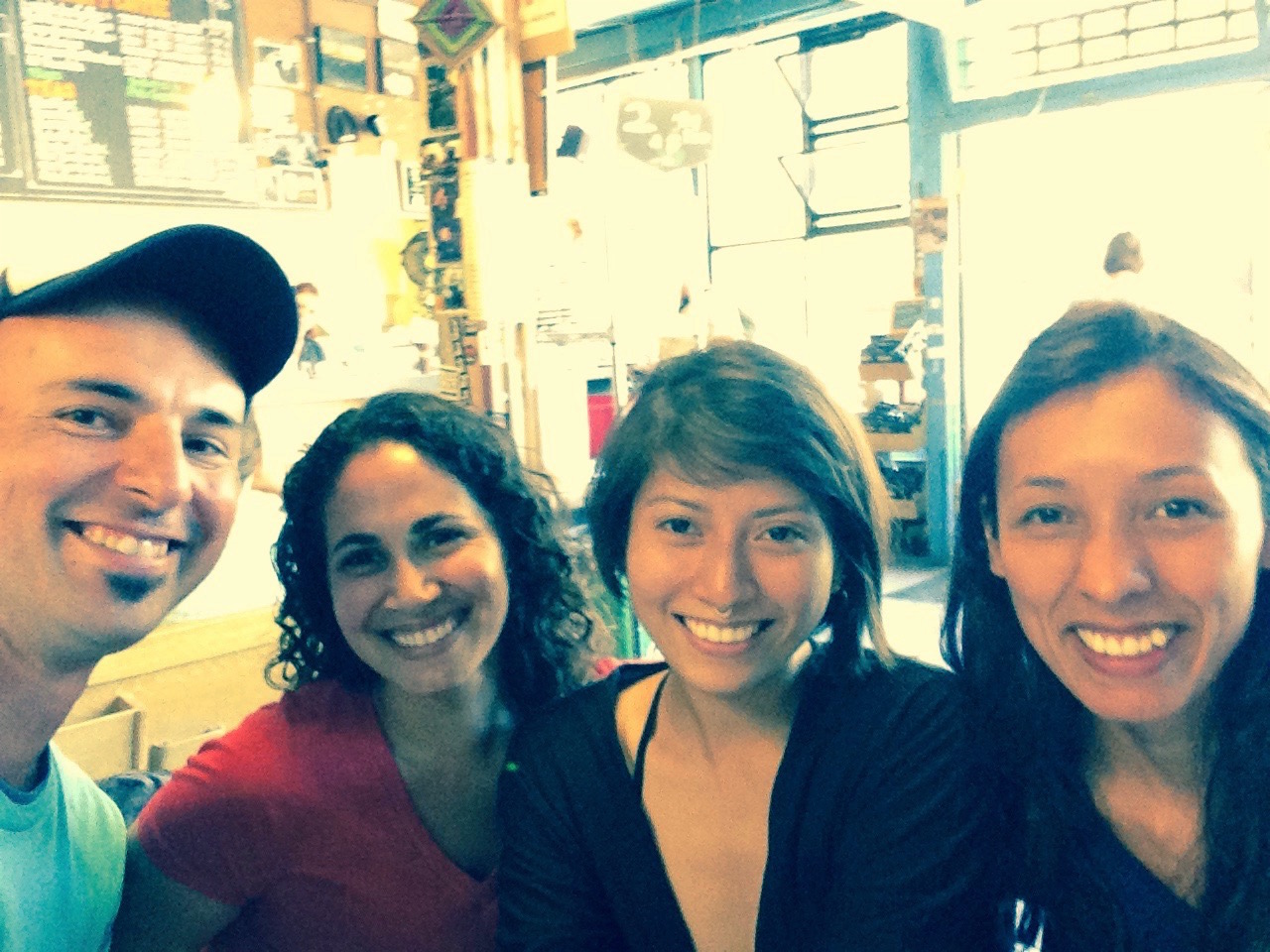} 
\caption{Celebrating Anastasia Chavez and Nicole Yamzon's new theorem on the Dehn--Sommerville relations \cite{CY}, Carolina Benedetti's visit to San Francisco, and Nicole's birthday.}
\end{figure}

\smallskip

\noindent 
\textsf{\textbf{Unconscious bias and discrimination.}} The mathematical society at large is relatively homogeneous, and students who do not fit neatly within its dominant cultures are often faced with unspoken but real obstacles. This is especially true for many underrepresented minorities, women and gender minorities, parents, disabled, low income, first-generation, and returning students among others.
On a professional level, 
these students face skepticism and their work is undervalued by professors and peers.\cite{CMA, LCMF, MDBGH, Y}
On a personal level, many feel pressured to leave their true selves at the door to try to fit in.


These are deep problems that we are far from solving. When they do arise, we talk about them openly, learn from them, and try to improve. We also emphasize that this is not just a problem for marginalized groups to address; every member of the mathematical society plays a role.

%
%



%
%
%

\smallskip

\noindent
\textbf{\textsf{The bigger picture.}} 
The SFSU-Colombia Combinatorics Initiative tries to instill an awareness of the tremendous power that mathematicians hold in today's society, and a collective belief in using our part of that power to do positive work. 
 Each individual gets to  decide what this means for them, whether or not they pursue careers in mathematics. For many of them, it means planting the seeds for the next generation of scientists in their communities.
Today's uneven representation in mathematics is largely, though not exclusively, a consequence of the uneven access to opportunities before students arrive to our college classrooms. Several of our students and alumni are doing powerful work with young people at various stages of their education. 
We are slowly integrating their initiatives through a network of mentorship that everyone contributes to and benefits from.


%


In Colombia, several participants in the initiative supervise undergraduate research projects and lead math  olympiad programs. Others partner with
the Clubes de Ciencia Colombia and the interactive science museum Parque Explora; thanks to them, ECCO 2016 will include a week-long workshop for students in Medell\'{\i}n public high schools, and a public math talk featuring the students' work which will be broadcast by television and streaming.

In the U.S., participants direct the San Francisco Math Circles, making great efforts to reach the (88\% ethnic minority) populations of the San Francisco Unified School District (SFUSD).
Through the GK-12 program, 50 Ph.D.-bound students spent 10 hours a week supporting the mathematics departments of various local public high schools. 
Other alumni teach mathematics full-time in local community colleges and high schools; the SFUSD claims that 70\% of its teachers have received training at SFSU.

In any initiative of this sort, one should share the decision making power and make oneself replaceable as soon as possible; 
I am inspired to be surrounded by a community of mathematicians who see themselves as agents of change in our scientific culture and in our societies.
They are doing extremely interesting mathematics, they are working hard to train the next generation of scientists -- a diverse, engaged, dynamic community that works to serve the needs of all -- and they are having a lot of fun in the process.

%
%
%

\begin{figure}[h]
\includegraphics[width=3.15in]{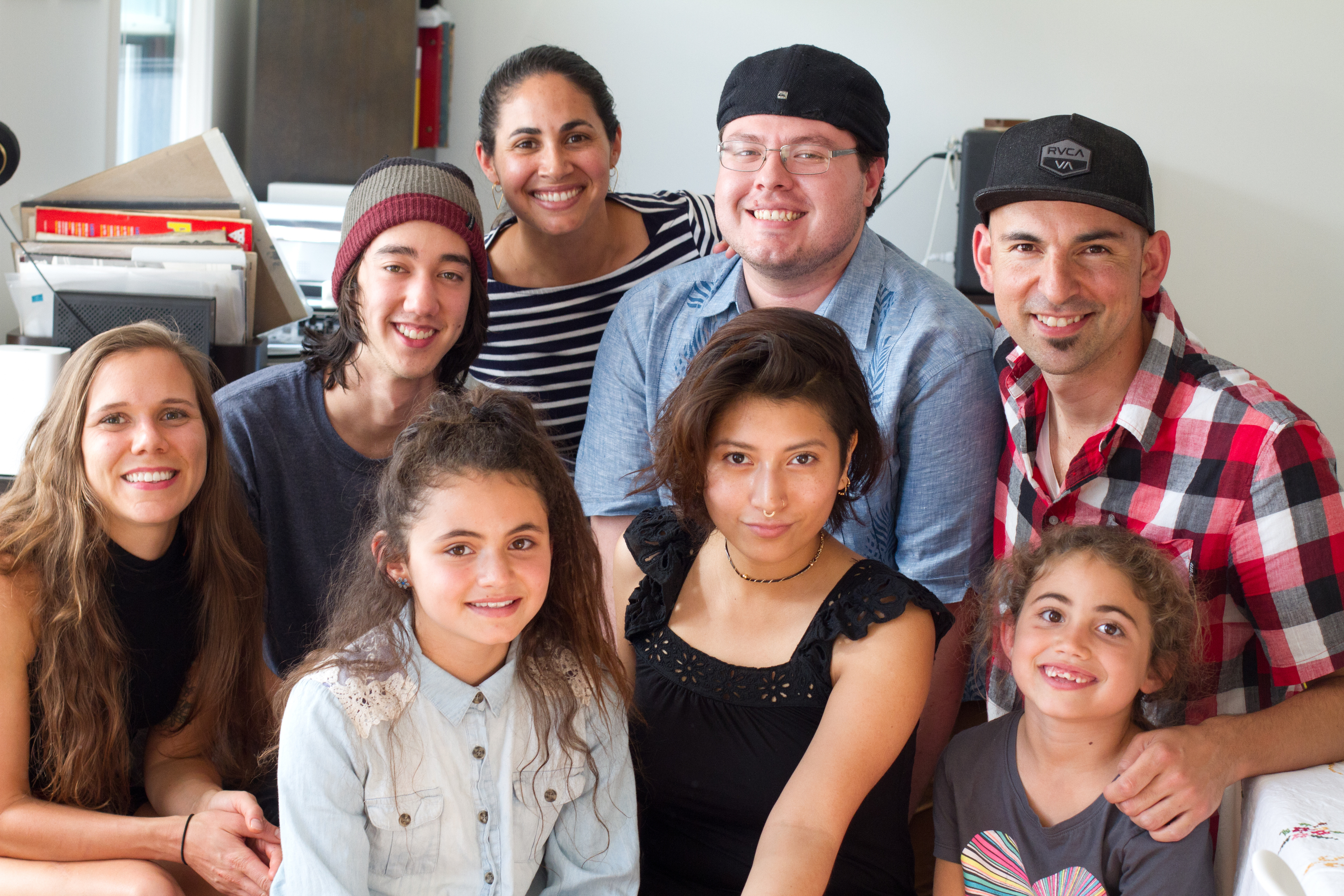} 
\end{figure}

\newpage

\noindent
\textbf{\textsf{More information.}}
Our webpage, containing additional information and resources, can be found at:

\begin{small}
\begin{center}
 \texttt{http://math.sfsu.edu/federico/sfsucolombia.html.}
\end{center}
\end{small}

\smallskip
\noindent
\textbf{\textsf{Acknowledgments.}}
I thank the AMS Committee on Education, AMS President Robert Bryant, and the editors of the Notices of the AMS  for inviting me to organize and share my ideas around a topic that is very challenging to me, partly because it is also very personal to me.

I am just one of many people doing this work at San Francisco State University and the in Colombia. I am indebted to my colleagues for their invaluable support of this project. The financial support by the NSF, the NIH, and CIMPA have also been instrumental.

If I have learned one thing, it is that I still have a lot to learn. I am extremely grateful to the teachers, colleagues, family members, friends, organizers, and futbolistas who have shaped this initiative; I must mention Natalia Ardila, Matthias Beck, Ben Braun, Dania Cabello, Jeff Duncan-Andrade, May-Li Khoe, Mar\'{\i}a de Losada, Amparo Mantilla de Ardila, Leticia M\'arquez-Maga\~na, Bob Moses, Ali Nesin, Gustavo Salazar, Kimberly Tanner, and Nicole Yamzon. Most importantly, I would like to thank my students, who teach me something new every day and give meaning to my mathematical career.

%

\scriptsize{

}

\begin{thebibliography}{99}

\bibitem{ABD}
Federico Ardila, Carolina Benedetti and Jeff Doker. Matroid polytopes and their volumes. \emph{Discrete and Computational Geometry} {\bf 43} (2010) 841-854. 

\bibitem{AC}
        Federico Ardila and Cesar Ceballos. 
Acyclic systems of permutations and fine mixed subdivisions of simplices. 
        \emph{Discrete and Computational Geometry.} {\bf 49} (2013) 485-510. 

\bibitem{AFR}
        Federico Ardila, Alex Fink and Felipe Rincon.
Valuations for matroid polytope subdivisions. 
        \emph{Canadian Journal of Mathematics} {\bf 62} (2010) 1228-1245.        



\bibitem{AR}
Federico Ardila and Amanda Ruiz.
When do two planted graphs have the same cotransversal matroid? 
\emph{Boletin de la Sociedad Matematica Mexicana} {\bf 16} (2010) 63-73. 

\bibitem{NSFreport}
Joan Burrelli. Academic institutions of minority faculty with science, engineering, and health doctorates. National Science Foundation InfoBrief, NSF 11-320, October 2011. \texttt{www.nsf.gov/statistics/infbrief/nsf11320/}

\bibitem{CBMSSurvey}
Richelle Blair, Ellen E. Kirkman, and James W. Maxwell. Statistical Abstract of Undergraduate Programs in the Mathematical Sciences in the United States. American Mathematical Society, 2013. 

\bibitem{CY}
Anastasia Chavez and Nicole Yamzon. The Dehn-Sommerville Relations and the Catalan Matroid. Preprint, 2015. \texttt{arXiv:1512.04513}

\bibitem{CMA}
Dolly Chugh, Katherine Milkman, and Modupe Akinola. Professors are prejudiced too. \emph{The New York Times Sunday Review}, pg. SR14, May 9, 2014.


\bibitem{CI}
Pauline Rose Clance and Suzanne Imes. The Imposter Phenomenon in High Achieving Women: Dynamics and Therapeutic Intervention. 
\emph{Psychotherapy: Theory, Research and Practice}. {\bf 15} (1978) 241-247

\bibitem{DF}
Harm Derksen and Alex Fink.
Valuative invariants for polymatroids.
\emph{Advances in Mathematics } {\bf 225} (2010), 1840--1892. 

\bibitem{D}
Carol Dweck. The secret to raising smart kids. \emph{Scientific American Mind} {\bf  18} (2007) 36-43.

\bibitem{G}
Daniel Z. Grunspan, Sarah L. Eddy, Sara E. Brownell, Benjamin L. Wiggins, Alison J. Crowe, and Steven M. Goodreau. 
Males Under-Estimate Academic Performance of Their Female Peers in Undergraduate Biology Classrooms.  PLOS ONE, {\bf 11} (2016), e0148405.

\bibitem{J}
Audrey Williams June. The invisible labor of minority professors. \emph{Chronicle of Higher Education}. Nov. 8, 2015.

\bibitem{LCMF}
Sarah-Jane Leslie, Andrei Cimpian, Meredith Meyer, Edward Freeland. Expectations of brilliance underlie gender distributions across academic disciplines. \emph{Science} {\bf 347} (2015) 22--265.

\bibitem{L}
Audre Lorde. Sister outsider: Essays and speeches. Crossing Press, 2012.
	



\bibitem{MLHA}
Joya Misra, Jennifer Hickes Lundquist, Elissa Holmes, and Stephanie Agiomavritis. The ivory ceiling of service work. \emph{Academe} {\bf 97} (2011): 22.


\bibitem{MOWA}
Kristen Monroe, Saba Ozyurt, Ted Wrigley, and Amy Alexander. 
Gender Equality in Academia: Bad News from the Trenches, and Some Possible Solutions. \emph{Perspectives on Politics} {\bf 6} (2008) 215-233

\bibitem{MDBGH}
Corinne A. Moss-Racusin, John F. Dovidio, Victoria L. Brescoll, Mark J. Graham, and Jo Handelsman. Science faculty's subtle gender biases favor male students
\emph{Proceedings of the National Academy of Science}  {109} (2012) 16474--16479.


\bibitem{Sz}
Martin Raussen and Christian Skau, Interview with Endre Szemerédi. Notices of the American Mathematical Society {\bf 60} (2013) 221 -- 231.

\bibitem{S}
Claude Steele. Whistling Vivaldi: And Other Clues To How Stereotypes Affect Us. New York : W.W. Norton and Company, 2010. 

\bibitem{Sc}
Martin A. Schwartz. The importance of stupidity in scientific research. \emph{Journal of Cell Science} {\bf 121} (2008) 1771.


\bibitem{T}
Kimberly Tanner. Structure matters: twenty-one teaching strategies to promote student engagement and cultivate classroom equity. \emph{CBE-Life Sciences Education} {\bf 12} (2013) 322-331.


\bibitem{Th}
William Thurston. On proof and progress in mathematics. \emph{Bulletin of the American Mathematical Society} {\textbf 30}   (1994) 161--177.

\bibitem{UScensus}
United States Census Bureau, 2014, accessed June 11, 2016. \texttt{https://www.census.gov/quickfacts/table/PST045215/00}


\bibitem{AMSSurvey}
William Yslas V\'elez, James Maxwell, and Colleen Rose. Report on the 2013-2014 new doctoral recipients. \emph{Notices of the American Mathematical Society}  {\bf 62} (2015) 771-781.

\bibitem{Y}
Ed Yong. XY Bias: How Male Biology Students See Their Female Peers. \emph{The Atlantic}

\bibitem{Z}
Estanislao Zuleta. Educaci\'on y democracia: un campo de combate. Corp. Tercer Milenio, Bogot\'a, 1995.

\end{thebibliography}
\end{document}